
 \magnification=1200
\input amssym.def
\input amssym.tex
 \font\newrm =cmr10 at 24pt
\def\bul{\raise .9pt\hbox{\newrm .\kern-.105em } }

 \def\fr{\frak}
 \font\sevenrm=cmr7 at 7pt

 \baselineskip 15pt
 
 \def\h{\hbox{ }}

 \def\r{{\fr r}}

 \def\n{{\fr n}}

 \def\ss{{\fr s}}
 
 \def\b{{\fr b}}
 \def\cc{{\fr c}}
 \def\hh{{\fr h}}

 \def\g{{\fr g}}

 \def\<{\le}
 \def\>{\ge}

 \def\s{{\h\subset\h}}
 
 \def\vs{\vskip }

 \def\mapright#1
  {\smash{\mathop
  {\longrightarrow}
  \limits^{#1}}}

 \def\kk#1{{\kern .4 em} #1}
 \def\vs{\vskip 1pc}

\hsize = 31pc
\vsize = 45pc
\overfullrule = 0pt
\rm

\overfullrule = 0pt
\rm

\centerline{\bf  Coadjoint structure of Borel subgroups and
their nilradicals} \vskip 1.5 pc \centerline{\it Bertram
Kostant}\vskip 1.5pc {\bf ABSTRACT}: Let
$G$ be a complex simply-connected semisimple Lie group and let $\g=
 \hbox{\rm Lie}\,G$. Let $\g = \n_- +\hh + \n$ be a triangular decomposition of
$\g$. One readily has that  $\hbox{\rm Cent}\,U(\n)$ is isomorphic to the ring 
 $S(\n)^{\n}$ of symmetric invariants. Using the cascade ${\cal B}$ of strongly orthogonal
roots, some time ago we  proved  that $S(\n)^{\n}$ is a polynomial
ring $\Bbb C[\xi_1,\ldots,\xi_m]$ where $m$ is the cardinality of ${\cal B}$. Using this  result
we establish  that the maximal  coadjoint of $N = \hbox{exp}\,\n$ has codimension $m$. 

Let $\b= \hh + \n$ so that the corresponding subgroup $B$ is a Borel subgroup of $G$. 
 Let $\ell = \hbox{rank}\,\g$.
Then in this paper we prove\vs {\bf Theorem.} {\it The maximal
coadjoint orbit of $B$ has codimension $\ell - m$ so that the
following statements are equivalent:
 $$\eqalign{&
(1)\,\,\,-1\,\,\hbox{is in the Weyl group of $G$ ( i.e., $\ell = m$)}\cr
&(2)\,\,\, B\,\, \hbox{has a nonempty open coadjoint orbit.} \cr}$$}

\indent {\bf Remark.}
A nilpotent or a semisimple group cannot have a nonempty open coadjoint orbit.
Celebrated examples where a solvable Lie group has a nonempty open coadjoint
orbit are due to Piatetski--Shapiro in his counterexample construction of a
bounded complex homogeneous domain which is not of Cartan type.

In a subsequent paper we will apply the results above to a
split real form $G_{\Bbb R}$. If
$\ell =1$  this gives rise (E. Stein) to the Hilbert transform. 

\noindent {\bf Keywords:}  coadjoint orbits, Borel subgroups, symplectic manifolds,
 root systems, structure theory

\noindent{\bf MSC 2010 codes:} 17B08, 17B22, 17B30, 22Exx, 53D05

\vskip 1.2pc
\centerline{\bf 1. Introduction}
\centerline{\bf {\bf Part A}: Recollection of some
of the relevant results of [K]}\vskip 1pc 

\noindent {\bf 1.1.} Let $G$ be a complex semisimple Lie group and let $\g = \hbox{Lie}\, G$. Let
$\g=
\n_-+\hh +
\n$ be a triangular decomposition of $\g$ and let $\ell$ be the rank of $\g$. Let
$\b =
 \hh
+ \n$ and $\b_- = \n_-+\hh$. Let $B,H,N$ be the subgroups of $G$ corresponding,
respectively, to $\b,\,\hh\,\n, $ so that $B$ is a Borel subgroup of $G$ and one has
the familiar semidirect product $B= H\ltimes N$.

Let ${\cal K}$ be the Killing form $(x,y)$ on $\g$. We may identify the dual $\b^*$
of $\b$ with $\b_-$, where $v\in \b^*$  and $z\in \b$; then $\langle v,z\rangle =
(v,z)$. With the obvious similar definitions the dual of $\n^*$ is identified
with $\n_-$. Let $\Phi_{\b}:\g\to \b_-$ be the projection defined by the
decomposition $\g = \b_- \oplus \n$, and let $\Phi_{\n}:\g\to \n_-$ be the
projection defined by the decomposition $\g= \n_-\oplus \b$. 
The coadjoint
representation $\hbox{Coad}_{\b}$ of $B$ may be given by 
$$\hbox{\rm Coad}_{\b}(g)v =\Phi_{\b} \hbox{\rm Ad}\,g\,\,v,$$ where $g\in B$ and $v\in \b_-$. For simplicity
we will write $g\bul\, v$ for
$\hbox{Coad}_{\b}(g)v$. The coadjoint representation
$\hbox{Coad}_{\n}$ of
$N$ is given  similarly where $\n$ replaces $\b$, $\n_-$ replaces $\b_-$, $N$
replaces $B$ and $\bul$ by $\cdot \,$. 

In order to prove the open $B$ coadjoint theorem mentioned in the abstract, we will
need the results of [K]. The results of [K] depend
upon definitions and properties of the cascade of orthogonal roots. These will be
used freely now, but for the convenience of the reader some of the definitions will
be here rccalled in Part A of the Introduction. 

The Killing form ${\cal K}$ induces a nonsingular bilinear form $(\mu,\nu)$ on $\hh^*$. 
 Let $\Delta\s \hh^*$ be the set of roots
corresponding to $(\hh,\g)$.  For each $\varphi\in \Delta$
let $e_{\varphi}\in \g$ be a corresponding root vector. The
root vectors can and will be chosen so that
$(e_{\varphi},e_{-\varphi}) = 1$ for all roots $\varphi$. If $\ss\s \g$ is any
 subspace stable under $ \hbox{\rm ad}\,\hh$ let
$$\Delta(\ss)=\{\varphi\in \Delta\mid e_{\varphi}\in \ss\}.$$ The set
$\Delta_+$ of positive roots is then chosen so that $\Delta_+ =
\Delta(\n)$, and one puts $\Delta_- = -\Delta_+$. 

Let ${\cal B}\s \Delta_+$ be the cascade of orthogonal
roots; See \S 1 in [K]. Then $\hbox{card}\,{\cal B} = m $ where $m$ is the maximal
number of  strongly orthogonal roots. See Theorem 1.8 in [K]. As in \S 2.1 in
[K], let $\r = \sum _{\beta\in {\cal B}}\Bbb C\,\,e_{\beta}$ so that $\r$ is an
$m$-dimensional commutative Lie subalgebra of $\n$. Let $R$ be the
$m$-dimensional commutative unipotent subgroup of $N$ corresponding to $\r$. As in
\S2.1 of [K], let
$\r_-\s \n_-$ be the span of $e_{-\beta}$ for $\beta\in {\cal B}$. For any $z\in
\r_-,\,\beta\in {\cal B}$, let $a_{\beta}(z)\in 
\Bbb C$
be defined so that $$z = \sum_{\beta\in
{\cal B}}\,a_{\beta}(z)\,e_{-\beta}\eqno (1.1)$$ and let $$\r_-^{\times} = 
\{\tau\in
\r_-\mid a_{\beta}(\tau)\neq 0,\,\,\forall
\beta\in {\cal B}\}.$$ As an algebraic subvariety of $\n_-$,
clearly 
$$\r_-^{\times} \cong (\Bbb C^{\times})^m.$$

For any $v\in \n_-$ let $O_v$ be the $N$-coadjoint orbit
containing $v$. Let $N_v\s N$ be the coadjoint isotropy subgroup at
$v$ and let $\n_v = \hbox{Lie}\,N_z$. Since the action is algebraic, 
$N_v$ is connected and as $N$-spaces $$O_v\cong N/N_v\,.  $$\vskip .5pc

The following result appears as Theorems 2.3 and 2.5 in [K]. \vs

 {\bf Theorem 1.1.} {\it Let $\tau\in \r_-^{\times} $. Then (independent
of $\tau$) $N_{\tau} = R$ so that (1.1) becomes $$O_{\tau}\cong
N/R.\eqno (1.2) $$
 Furthermore if
$\tau,\tau'\in \r_-^{\times} $ are distinct,  then $O_{\tau}\cap O_{\tau'}=
\emptyset$ so that one has a disjoint union
$$N\cdot \r_-^{\times}  = \cup_{\tau\in \r_-^{\times} }\,O_{\tau}.\eqno
(1.3)$$}\vskip .5pc
 Since $B$ normalizes $N$ there is  a natural action on $B$ on 
$\n_-$. We refer to this
as the $\n_-$ action of $B$. Explicitly since the latter extends $\hbox{Coad}_{\n}$ we
denote the $\n_-$ action of $B$ of $b\in B$ on $v\in \n_-$ by $b\cdot v$, and note
that $$b\cdot v=
\Phi_{\n}\circ \hbox{Ad}\,b\,\,v.\eqno (1.4)$$

The Zariski open subvariety $\r_-^{\times} \s \r_-$ is stable
 under the $\n_-$ action of $H\s B$
and in fact $H$ operates transitively on $\r_-^{\times} $ so that $\r_-^{\times} $ is
isomorphic to a homogeneous space for $H$.

 \vskip 6pt The action
of $H$ on $\r_-^{\times} $ extends to an action of $H$ on the corresponding set
$\{O_{\tau},\,\tau\in \r_-^{\times} \}$ of $N$-coadjoint orbits. Since $H$
normalizes $N$ the following statement is obvious. See \S 2.2 in [K].

 \vskip 6pt {\bf
Proposition 1.2.} {\it For any
$\tau\in \r_-^{\times} $ and
$a\in H$ one has
$$O_{{\hbox{\sevenrm Ad}}\,a(\tau)} = a\cdot O_{\tau}. $$}
\indent Let
$$X=\cup_{\tau\in \r_-^{\times} }\,O_{\tau}\eqno (1.5)$$ so that, by 
(1.3), 
the union (1.5) is disjoint. Furthermore since $B=NH$  we note that $X$ is an
orbit of $\n_-$ action of $B$ and hence is Zariski open and dense in its closure.  
The following
is one of the main theorems in [K]. See Theorem 2.8. in [K]. \vs {\bf Theorem
1.3.} {\it $X$ is Zariski dense in $\n_-$ so that $$\overline {X} = \n_-.\eqno
(1.6)$$}

  Let $W$ be the Weyl group of $\g$ operating in $\hh$
and
$\hh^*$. Reluctantly submitting to common usage, let $w_o$ be
the long element of $W$.  For any $\beta\in {\cal B}$ let $s_{\beta}\in W$ be the
reflection defined by $\beta$. Since the elements in ${\cal B}$ are orthogonal
to one another, the reflections $s_{\beta}$ evidently commute
with one another. The long element
$w_o$ of the Weyl group
$W$ is given in terms of the product of these commuting
reflections. In fact one has (see Proposition 1.10 in [K]). \vs {\bf
Proposition 1.4.} {\it One has
$$w_o= \prod_{\beta\in {\cal B}}\,\,s_{\beta}, \eqno (1.7)$$ noting that the
order of the product is immaterial because of commutativity.} 

\vs {\bf Remark.} Note that it
is immediate from Proposition 1.4 that $-1\in W$ if and only if $m=\ell$.\vskip 1.5pc  
\centerline{\it {\bf Part B}. Statement of new results to be proved in \S2.}\vskip
1pc 
{\bf 1.2}  With respect to the coadjoint structure of $N$ one has \vs {\bf Theorem
1.5.} {\it A maximal coadjoint orbit of $N$ has codimension $m$ in $\n^*$ where $m$ is the maximal
number of strongly orthogonal roots. In particular it has codimension $\ell$ (the rank of $\g$) in
$\n^*$ if and only if $-1\in W$.}

\vs We now deal with the coadjoint representation $\hbox{Coad}_{\b}$ of
the Borel subgroup $B$. As noted in \S1.1,  we have identified $\b^*$ with $\b_-$. For any
$w\in
\b_-$ let
$B_w$ be the $B$-coadjoint isotropy group of $B$ at $w$ and let $\b_w = \hbox{Lie}\,B_w$. We refer
to
$\b_w$ as the $B$-coadjoint isotropy algebra at $w$. Let
${\cal O}_w\s
\b_-$  be the
$B$-coadjoint orbit of
$w$ so that $${\cal O}_w \cong B/B_w. $$ 

 Of course $\hh \s \b$. Let $\hh^o\s \hh$ be
the orthogonal subspace to $\Bbb C {\cal B}\s \hh^*$ in $\hh$. Then $$\hbox{\rm dim}\,\hh^o = 
\ell - m.\eqno (1.8)$$ \vskip .5pc Now let $\Phi_{\hh}:\g \to \hh$ be the projection defined by the
triangular decomposition $\g = \n_- +\hh +\n$. Obviously any $w\in \b_-$ can be uniquely written
 as
$$w = v + x,\eqno (1.9)$$ where $v=\Phi_{\n}w$ and $x= \Phi_{\hh}w$. It is immediate from Theorem
1.3 that $\b_-'$ is open in $\b_-$ where $\b_-'$ is defined by $$\b_-' = \{w\in \b_-\mid
\Phi_{\n}w\in X\}.\eqno (1.10)$$ We will prove \vs{\bf Theorem 1.6.} {\it Let $w\in \r_-^{\times}$.
Then the $B$-coadjoint isotropy algebra of $\b$ at $w$ is
$\hh_{\cal}^o$. Furthermore ${\cal O}_w$ is a maximal $B$-coadjoint orbit so that $\ell- m$ is the
codimension of any maximal $B$-coadjoint orbit. Furthermore
${\cal O}_y$ is a maximal $B$-coadjoint orbit of $B$ for any $y\in \b_-'$.}\vs As a consequence we
now have \vs {\bf Theorem 1.7.} {\it A Borel subgroup of a complex
semisimple Lie group has a nonempty (and hence unique) open coadjoint orbit if and only if $m=\ell$,
that is if and only if
$-1$ is in the Weyl group $W$. Furthermore, explicitly, the unique nonempty open coadjoint orbit
of $B$ is
$\b_-'$.}\vskip 1.5pc \centerline{\bf 2. Proofs of stated results and additional results}\vskip
1pc {\bf 2.1.} The proof of  Theorem 1.5 depends on the density of $X$ and standard arguments
about dimensions of orbits.
\vs {\bf Proof of Theorem 1.5.}  Let $v\in X$. Since $\hbox{dim}\, R = m$ it follows from (1.2) and (1.5)
that the codimension of $O_v$ in $\n_-$ is $m$. But if $w\in \n_-$,  then by (1.6) there exists a
sequence
$v_k\in X$ which converges to $w$. But then by the compactness of subspaces of dimension $m$ there
exists a subsequence $v_{k'}$ such that $\n_{v_{k'}}$ converges to an $m$-dimensional subspace of
$\n_w$. Hence the codimension of $O_w$ in $\n_-$ is greater than or equal to  $m$. \hfill QED

\vskip .5pc
Recall (1.4). \vskip 6pt {\bf Theorem 2.1.} {\it Let $w\in \r_-^{\times}$. Then the isotropy subalgebra
of
$\b$ at $w$ with respect to the $\n_-$ action of $B$ is $\hh^o + \r$.}

\vskip 6pt {\bf Proof.} Let
$\cc_{w}$ be the isotropy subalgebra of
$\b$ at $w$ with respect to the $\n_-$ action of $B$. Since $B= NH$ and since the $\n_-$ action
of $B$ extends $\hbox{Coad}_{\n}$ clearly $(\hh^o + \r) \s \cc_{w}$. But $$B\cdot w = X\eqno (2.1)$$ by
Proposition 1.2 and Theorem 1.1. But $\hbox{dim}\,(\hh^o + \r) = \ell$ since clearly 
$$\hbox{dim}\,\hh^o = \ell-m\eqno (2.2)$$
and
$\hbox{dim}\,\r= m$. But
$\hbox{dim}\, X = \hbox{dim}\,\n$ by (1.6). Thus $\hbox{dim}\,\cc_{w} = \ell$. This implies $\hh^o + \r =
\cc_{w}$. \hfill QED\vskip 4pt
Comparing $\hbox{Coad}_{\b}$ with the $\n_-$ action of $B$, we note that if $w\in
\n_-$ and $b\in B$ (recall the notation of \S 1.1.), one has $$b\bul\, w = b\cdot w + \Phi_{\hh}\,
 b\bul\, w\eqno (2.3)$$ This is immediate since $\Phi_{\b} = \Phi_{\n} + \Phi_{\hh}$.  This latter sum
also yields the following infinitesimal analogue of (2.3). Let $x\in \b$ and $w\in \n_-$;  then
$$x\bul\, w = x\cdot w + \Phi_{\hh} x\bul \, w \eqno (2.4)$$ where $x\bul\, w = \Phi_{\b}\,[x,w]$ and
$x\cdot w =\Phi_{\n}[x,w]$. But now (2.1) and (2.3) imply

\vs {\bf Proposition 2.2.} {\it Let
$w\in
\r_-^{\times}$. Then
$$\Phi_{\n}\,B\bul\, w =  X.\eqno (2.5)$$}\vskip .5pc We also note that the fixed point set $\b_-^B$
of the coadjoint action of $B$ is immediately given by $$\b_- ^B = \hh.\eqno (2.6)$$ 
One now has \vs
{\bf Theorem 2.3.} {\it Let
$w\in
\r_-^{\times}$. Then the isotropy algebra, $\b_w$, of  $\b$ at $w$, with respect to the coadjoint
action of
$B$ is $\hh^o$ so that, by (2.2), $\ell- m$ is the codimension of the coadjoint orbit $B\bul\, w$.} 
\vs
{\bf Proof.} Let $x\in \b_-$. Then, by Theorem 2.1 and (2.4) one has $x\in \b_w$ if and only if 
(1) $x\in 
\hh^o + \r$ and (2) $\Phi_{\hh} [x,w]=0$. But clearly $\hh^o\s \b_w$. Thus to prove $\b_w=\hh^o$
it suffices to show that if $$0\neq x\in\r\,\,\,  \hbox{then}\,\,\Phi_{\hh} [x,w]\neq 0 .\eqno
(2.7)$$ But since the all the  coefficients $a_{\beta}(w)$ in (1.1) are not zero and the set
$\{\beta^{\vee}\in
\hh
\mid
\beta\in {\cal B}\}$ are linearly independent, one has (2.7). \hfill  QED\vs {\bf 2.2.} As an immediate
consequence of (2.6) one has 

\vs {\bf Proposition 2.4.} {\it Let $w\in \b_-$ and $z\in \hh$; then for
any
$v\in {\cal O}_w$ one has $v+z\in {\cal O}_{w+ z}$ and the map 
$${\cal O}_w\to {\cal O}_{w+z},\,\,\,\,v \mapsto v+z\eqno (2.8)$$ is a $B$-isomorphism of
$B$-coadjoint orbits.} 

\vs As an immediate
consequence of Proposition 2.4, one has \vs

  {\bf Proposition 2.5.} {\it Let
$w\in \b_-$ and $z\in \hh$. Then $$\b_{w} = \b_{w+z}\eqno (2.10)$$ so that obviously in particular
$$\hbox{\rm dim}\,\b_w=\hbox{\rm dim}\,b_{w+z}.\eqno (2.10)$$}
\indent  From the open density of $X$ in $\n_-$, it is obvious,
recalling (1.10), that $\b_-'$ is Zariski open and dense in $\b_-$. \vs {\bf Theorem 2.6} {\it
For any $y\in \b_-' $ the dimension of $\b_y$ in $\b$ is $\ell-m$ where we recall that $m$ is the
cardinality of a maximal set of strongly orthogonal roots. Furthermore ${\cal O}_y$ is a maximal
$B$-coadjoint orbit and has codimension $\ell-m$ in $\b_-$. That is,  the codimension of 
any $B$-coadjoint orbit ${\cal O}$ is greater than or equal to $\ell-m$.}

\vs  {\bf Proof.} Let $v =
\Phi_{\n} y$ so that $v\in X$. But then by (2.1) and (2.3) there exists $w\in \r_-^{\times}$
and $b\in B$ such that $\Phi_{\n}b\bul \, w = v$. But by conjugacy and Theorem 2.3, $
\hbox{dim}\,\b_{b\bul \, w} = \ell -m$. But if $z = y- b\bul \, w$, then $z\in \hh$. But then $
\hbox{dim}\,\b_y =
\ell-m$ by Proposition 2.5. But then using the limit argument in the proof of Theorem 1.5, one has
$\hbox{dim}\,\b_x \geq \ell-m$ for any $x\in \b_-$. \hfill QED

\vs Applying Theorem 2.6 to the case where
$\ell=m$ we come to the statement and proof of our main theorem. \vs {\bf Theorem 2.7.} {\it Let
$G$ be a complex semisimple Lie group and let $B$ be a Borel subgroup of $G$. Then $B$ has a
(necessarily unique) nonzero open coadjoint orbit, ${\cal O}$, if and only if $-1$ is in the Weyl
group of $G$, or equivalently, there exists $\ell$ orthogonal roots where $\ell$ is the rank of
$G$. Furthermore in the notation above $${\cal O} = \b_-'\eqno (2.11)$$ where $\b_-'$ is given by
(1.10). In particular $\r_-^{\times }\s {\cal O}$. (See \S 1.1.)} 

\vs {\bf Proof.} By Theorem 2.6
$B$ has an open coadjoint orbit ${\cal O}$ in $\b_-$ if and only if $m=\ell$. The orbit is
necessarily unique since the action of $B$ is  linear algebraic and $\b_-$ is affine
irreducible. But now $\b_-'\s {\cal O}$ by Theorem 2.6. On the other hand,  ${\cal O}\s \b_-'$ by
Proposition 2.2. This proves (2.11). \hfill QED

\vskip 4pt {\bf Remark 2.8.} If $G$ is simple, then $B$ has a
nonempty open coadjoint orbit if and only if $\g$ is of type $B_{\ell},\, C_{\ell},\,
D_{\ell}\,\,\,
\hbox{with $\ell$-even},\,\,G_2,\,F_4,\,E_7,\,E_8$. Except for the case $D_{\ell},\,\,\,
\hbox{with $\ell$-even}$, the other cases follow from the fact that there exists no 
outer automorphism. This fact readily implies $-1$ is in the Weyl group. \vskip
6pt
\centerline{\bf References} \vskip4pt

\item {[K]} B. Kostant, The cascade of orthogonal roots and the
coadjoint structure of the nilradical of a Borel subgroup of a
semisimple Lie group, to appear in the {\it  Moscow Math Journal},
edited by Victor Ginzburg, in memory of I. M. Gelfand, Spring 2012.
Also in arXiv:1101.5382 [math RT] \vskip 2pt

\noindent Bertram Kostant (Emeritus)

\noindent Department of Mathematics

\noindent Cambridge, MA 02139

\noindent kostant@math.mit.edu

\end

\end